\documentstyle[12pt,amscd,amssymb]{amsart}

\def\sw#1{{\sb{(#1)}}}
\def\su#1{{\sp{(#1)}}} 
\def\tens{\mathop{\otimes}}

\def\endproof{\hbox{$\sqcup$}\llap{\hbox{$\sqcap$}}\medskip}
\def\<{{\langle}}
\def\>{{\rangle}}
\def\ra{{\triangleleft}}
\def\la{{\triangleright}} 
\def\eps{\epsilon}

\def\note#1{{}}
\def\can{{\rm can}}
\def\note#1{}
\def\M{{\bf M}}
\def\CF{{\cal F}}
\def\ev{{\rm ev}}
\def\coev{{\rm coev}}

\def\Label{\label}
\begin{document}

\newtheorem{proposition}{Proposition}[section]
\newtheorem{lemma}[proposition]{Lemma}
\newtheorem{corollary}[proposition]{Corollary}
\newtheorem{theorem}[proposition]{Theorem}

\theoremstyle{definition}
\newtheorem{definition}[proposition]{Definition}
\newtheorem{example}[proposition]{Example}

\theoremstyle{remark}
\newtheorem{remark}[proposition]{Remark}

\newcommand{\Section}{\setcounter{definition}{0}\section}
\renewcommand{\theequation}{\thesection.\arabic{equation}}
\newcounter{c}
\renewcommand{\[}{\setcounter{c}{1}$$}
\newcommand{\etyk}[1]{\vspace{-7.4mm}$$\begin{equation}\Label{#1}
\addtocounter{c}{1}}
\renewcommand{\]}{\ifnum \value{c}=1 $$\else \end{equation}\fi}

\title[Properties of Entwined Modules]{Frobenius 
Properties and Maschke-type Theorems for
 Entwined Modules}
\author{Tomasz Brzezi\'nski}
\address{Department of Mathematics, University of York, Heslington, 
York YO1 5DD, U.K.}
\email{tb10@@york.ac.uk}
\urladdr{http//www.york.ac.uk/\~{}tb10}
\thanks{Lloyd's of London Tercentenary Fellow.}
\thanks{On leave from: Department of Theoretical Physics, University of
\L\'od\'z, Pomorska 149/153, 90-236, \L\'od\'z, Poland}
\date{June 1998}
\begin{abstract}
Entwined modules arose from the coalgebra-Galois theory. They are a 
generalisation of unified Doi-Hopf modules. In this paper,  Frobenius
properties and Maschke-type theorems, known for Doi-Hopf modules 
are extended to the case of entwined
modules. 
\end{abstract}
\maketitle
\section{Introduction}
The aim of this article is to show that results of recent papers
\cite{CaeMil:Doi} and \cite{CaeMil:Mas} concerning  
Frobenius properties and a Maschke-type theorem for Doi-Hopf modules 
\cite{Doi:uni}
\cite{Kop:var} hold for the more
general class of modules, known as {\em entwined modules}
\cite{Brz:mod}. 
These are
modules of an algebra and comodules of a coalgebra such that the action
and the coaction satisfy certain compatibility condition. Unlike
Doi-Hopf modules, entwined modules are defined purely using the
properties of an algebra and a coalgebra combined into an {\em entwining
structure}. There is no need for a `background' bialgebra, which is
an indispensable part of Doi-Hopf construction. 
The bialgebra-free
formulation - apart from being more general - has also remarkable
self-duality property, which essentially implies that for every statement
involving module structure of an entwined module there is a corresponding
statement involving its comodule structure. 
 As an
illustration of this phenomenon we give two Maschke-type theorems for
entwined modules. 

The paper is organised as follows. In Section~2 we recall definitions
 and give examples of entwining structures and 
entwined modules. In Section~3 we introduce integrals for entwining
structures and analyse Frobenius properties of
entwined modules generalising the results of \cite{CaeMil:Doi}. 
Finally, in Section~4 we state a Maschke-type theorem
for entwined modules and derive its dual form.

We work over a commutative ring $k$. We  assume that all
the algebras are over $k$ and  unital, and the coalgebras are over $k$
and  counital. Unadorned tensor product is over $k$. For any $k$-modules
$V, W$ the symbol  ${\rm Hom}(V,W)$ denotes the $k$-module
of $k$-module maps $V\to W$, the identity map $V\to V$ is
 denoted by $V$, $V^* = {\rm Hom}(V,k)$, and $\ev_V :V\otimes V^*\to k$
denotes the evaluation map, i.e. $\ev_V :v\otimes v^* 
\mapsto \<v,v^*\> = v^*(v)$. We also implicitly identify $V$ 
with $V\otimes k$
and $k\otimes V$ via the canonical isomorphisms.

For a $k$-algebra $A$ we use $\mu$ to denote the product as a map and 
$1_A$ to denote  unit both as an element of $A$ and as a map $k\to
A$, $\alpha\to \alpha 1_A$. $\M_A$ (resp.\ ${}_A\M$) denotes 
the category of right (resp.\ left) $A$-modules. The morphisms in this
category are denoted by ${\rm Hom}_A(M,N)$ (resp.\ ${}_A{\rm Hom}(M,N)$).
For any $M\in \M_A$ (resp.\ $M\in {}_A\M$), the symbol $\rho_M$ (resp.\
${}_M\rho$) denotes the action as a map (on elements the action  is
denoted by a dot). We often write $M_A$ (resp.\ ${}_AM$) to
indicate in which context the $A$-module $M$ appears.
If $A,B$ are $k$-algebras and $M,N$ are $(A,B)$-bimodules then ${}_A{\rm
Hom}_B(M,N)$ denotes the set of $(A,B)$-bimodule maps.

For a
$k$-coalgebra $C$ we use $\Delta$ to denote the coproduct and $\eps$ to
denote the counit. Notation for comodules is similar to that for modules
but with subscripts replaced by superscripts, i.e. $\M^C$ is the
category of right $C$-comodules, $\rho^M$ is a right coaction etc. We
use the Sweedler notation for coproducts and coactions, i.e. $\Delta(c)
= c\sw 1\otimes c\sw 1$, $\rho^M(m) = m\sw 0\otimes m\sw 1$ (summation
understood). The symbol $\la$ (resp.\ $\ra$) stands for the standard left
(resp.\ right) 
action of the (convolution product) algebra
$C^*$ on $C$ induced by the coproduct in $C$.

\section{Preliminaries on entwining structures}
\begin{definition}
\Label{ent}
An {\em  entwining
structure} (over $k$) is a 
triple $(A,C)_\psi$ consisting of a $k$-algebra $A$, a $k$-coalgebra
$C$ and a $k$-module map $\psi: C\tens A\to A\tens C$ satisfying
$$
\psi\circ(C\tens \mu) = (\mu\tens C)\circ (A\tens\psi)\circ(\psi\tens A),
\qquad \psi\circ (C\tens 1_A) = 1_A\tens C,
$$
$$
(A\tens\Delta)\circ\psi = (\psi\tens
C)\circ(C\tens\psi)\circ(\Delta\tens A), \qquad (A\tens \eps)\circ\psi =
\eps\tens A.
$$
A morphism of entwining structures is a pair
 $(f,g):(A,C)_\psi\to
(\tilde{A}, \tilde{C})_{\tilde{\psi}}$, where $f:A\to \tilde{A}$ is
an algebra map, $g: C\to \tilde{C}$ is a coalgebra map, and
$(f\otimes g)\circ\psi =
\tilde{\psi}\circ(g\otimes f)$. 
\end{definition}

For $(A,C)_\psi$ we  use the
notation $\psi(c\otimes a) = a_\alpha \otimes c^\alpha$ (summation
over a Greek index understood), for all $a\in
A$, $c\in C$.
The notion of an entwining structure was  introduced in
\cite[Definition~2.1]{BrzMa:coa}. It is self-dual in the sense that
conditions in Definition~\ref{ent} are invariant under the operation
consisting of  interchanging
of $A$ with $C$, $\mu$ with $\Delta$, and  $1_A$ with $\eps$, and
reversing the order of maps. Below are two main classes of examples
of entwining  structures. 

\begin{example}
Let $H$ be a bialgebra, $C$ a right
$H$-module coalgebra, and $A$ a right $H$-comodule algebra. Then 
$C$ and $A$ are
entwined by $\psi: C\otimes A\to 
A\otimes C$, $c\otimes a\mapsto a\sw 0\otimes  c\cdot a\sw 1$. The
corresponding entwining structure $(A,C)_\psi$ is called an {\em
entwining structure associated to a Doi-Hopf
datum $(A,C,H)$}. 
\Label{Doi-Hopf}
\end{example}
\begin{example}[\cite{BrzHaj:coa}]
Let $C$ be a coalgebra, $A$ an algebra and a right $C$-comodule. Let 
$B:= \{b\in A\; | \; \rho^A(ba) =
b\rho^A(a)\}$ and assume that the canonical left $A$-module, right
$C$-comodule map
$
can:A\otimes _BA\to A\otimes C$, $a\otimes a'\mapsto a\rho^A(a')$,
is bijective. Let $\psi:C\otimes A\to A\otimes C$ be a $k$-linear map
given by
$
\psi(c\otimes a) = can(can^{-1}(1_A\otimes c)a).
$
Then $(A,C)_\psi $ is an entwining structure. The 
extension $B\subset A$ is called a
{\em coalgebra-Galois extension} (or a {\em $C$-Galois extension})  and is
denoted by $A(B)^C$. $(A,C)_\psi$ is the
{\em canonical entwining structure} associated to $A(B)^C$.
\label{can.ex}
\end{example} 

\begin{definition}
Let $(A,C)_\psi$ be an entwining structure. An (entwined) 
{\em $(A,C)_\psi$-module}
is a right $A$-module, right $C$-comodule $M$ 
such that 
$$
\rho^M\circ\rho_M = (\rho_M\tens C)\circ(M\tens\psi)\circ(\rho ^M\otimes
A), 
$$
(explicitly: 
$
\rho^M(m\cdot a) = m\sw 0\cdot a_\alpha\tens m\sw 1^\alpha$, $ \forall
a\in A, m\in M$).
A morphism of $(A,C)_\psi$-modules is a right $A$-module map
which is also a right $C$-comodule map. The category of
$(A,C)_\psi$-modules is denoted by $\M_A^C(\psi)$.
\label{def.mpsi}
\end{definition}
The category $\M_A^C(\psi)$ was introduced and
studied in \cite{Brz:mod}. Modules associated to the entwining
structure in Example~\ref{Doi-Hopf} are unifying Hopf modules or
 Doi-Hopf modules introduced in \cite{Doi:uni}, \cite{Kop:var}. On the
other hand, 
entwined modules associated to the entwining structure in
Example~\ref{can.ex} do not seem to be of the Doi-Hopf type.
The following  example  is a special case of the construction in 
\cite[Section~3]{Brz:mod}.

\begin{example} Let $(A,C)_\psi$ be an entwining structure. Then 

(1) If $M$ is a right $A$-module then $M\otimes C\in \M^C_A(\psi)$ with
the coaction $M\otimes \Delta$ and the action $(m\otimes c)\cdot a =
m\cdot\psi(c\otimes a)$, for all $a\in A, c\in C$ and $m\in M$. In
particular $A\otimes C\in \M^C_A(\psi)$.

(2) If $V$ is a right $C$-comodule then $V\otimes A\in \M^C_A(\psi)$
with the action $V\otimes\mu$ and the coaction $v\otimes a\mapsto v\sw
0\otimes \psi(v\sw 1\otimes a)$ for any $a\in A$ and $v\in V$. In
particular $C\otimes A\in \M^C_A(\psi)$.
\label{ex.psi.modules}
\end{example}

\section{Integrals and Frobenius properties of entwined modules}
An entwining structure $(A,C)_\psi$ is said to be {\em factorisable} if
there exists a unique map $\bar{\psi} :A\otimes C^*\to C^*\otimes A$ 
such that the following diagram
$$
\begin{CD}
C\otimes A\otimes C^* @>{C\otimes\bar{\psi}}>> C\otimes C^*\otimes A\\
@VV{\psi\otimes C^*}V                          @VV{\ev_C\otimes A}V\\
A\otimes C\otimes C^*  @>{A\otimes \ev_C}>>  A
\end{CD}
$$
commutes. For any $a\in A$ and $\xi\in C^*$
we write $\bar{\psi}(a\otimes\xi) =
\xi_i\otimes a^i$ (summation over repeated index understood). 
For example, $(A,C)_\psi$ associated to a Doi-Hopf datum in
Example~\ref{Doi-Hopf} is factorisable, provided $C$ is a projective
$k$-module. Also, any $(A,C)_\psi$ with $C$ being finitely generated
projective $k$-module is factorisable. 

Let $(A,C)_\psi$ be a
factorisable entwining structure and 
let $B=C^{*op}$, i.e. $\<c,bb'\> =
\<c\sw 2,b\>\<c\sw 1,b'\>$, for all $b,b'\in B$, $c\in C$.
As explained in \cite[Proposition~2.7]{BrzMa:coa}, the triple
$(A,B,\bar{\psi})$ defines a factorisation structure (cf.
\cite[Equations~(7.10)]{Ma:book}). This allows one to construct a 
generalised {\em smash
product} (factorised) algebra $B\#_{\bar{\psi}} A$ on $B\otimes A$ with the product $(b\otimes
a)(b'\otimes a') = b\bar{\psi}(a\otimes b')a'$, for all $a,a'\in A$,
$b,b'\in B$. The maps $A\hookrightarrow B\#_{\bar{\psi}} A$, $a\mapsto
1_B\otimes a$ and $B\hookrightarrow B\#_{\bar{\psi}} A$, $b\mapsto
b\otimes 1_A$ are algebra inclusions. Every $B\#_{\bar{\psi}} A$-module
is viewed as an $A$ or $B$ module via these maps. 
Given such $B\#_{\bar{\psi}} A$ we
equip the $k$-module  ${\rm Hom}(B,A)$ with the structure
of an $(A,B\#_{\bar{\psi}} A)$-bimodule via 
$(a\cdot f\cdot(b\otimes a'))(b') =
af(bb'_i){a'}^i$, for all $a,a'\in A$, $b,b'\in B$, $f\in {\rm
Hom}(B,A)$. 
A $k$-module map $\lambda :
B\to A$ such that for all $a\in A$, $a\cdot\lambda = \lambda\cdot a$, is
called an {\em integral} in $B\#_{\bar{\psi}} A$.
The $k$-module of integrals in $B\#_{\bar{\psi}}A$ is denoted by ${\rm Int}(B\#_{\bar\psi}A)$.

Finally, we recall the notion of a Frobenius extension introduced in
\cite{Kas:pro} and \cite{NakTsu:fro}

\begin{definition}
Let $X$ be an algebra and $A$ its subalgebra. The extension $A\subset X$
is called a {\em Frobenius extension} (of the first kind) iff $X$ is a
finitely generated projective right $A$-module and $X \cong 
{\rm Hom}_A(X,A)$ as $(A,X)$-bimodules. The $(A,X)$-bimodule
structure of ${\rm Hom}_A(X,A)$ is given by $(a\cdot f\cdot x)(x') =
af(xx')$, for all $a\in A, x,x'\in X$ and $f\in {\rm Hom}_A(X,A)$.
\end{definition}
\begin{proposition}
Let $(A,C)_\psi$ be an entwining structure and let $B = C^{*op}$. Assume that  $B$ is
a  finitely generated projective $k$-module and $A$ is a faithfully flat
$k$-module, and let $X= B\#_{\bar{\psi}} A$. Then the following are 
equivalent:

(1) The functor ${\rm Hom}(B,-):\M_A\to \M_X$ is the left adjoint
of the  functor $\CF: \M_X\to \M_A$ induced by $A\hookrightarrow X$.

(2) The extension $A\subset X$ is Frobenius.

(3) $X$ is isomorphic to ${\rm Hom}(B,A)$ as an $(A,X)$-bimodule.

(4) There exists an integral $\lambda\in {\rm Int}(B\#_{\bar\psi}A)$ such that the
map $\phi: B\otimes A\to {\rm Hom}(B,A)$, $b\otimes a\mapsto
\lambda\cdot (b\otimes a)$ is bijective.
\label{prop.frobenius1}
\end{proposition}
\begin{proof}
(1) $\Leftrightarrow$ (2). ${\rm Hom}(B,-)$
is the right adjoint of $\CF$ since for all $M\in \M_X$, $N\in
\M_A$ there is a natural isomorphism $\eta_{M,N}: {\rm Hom}_X(M, {\rm
Hom}(B,N)) \to {\rm Hom}_A(M,N)$, $\eta_{M,N}(f)(m) =
f(m)(1_B)$, for all $m\in M$. Its inverse is 
$\eta^{-1}_{M,N}(g)(m)(b) = g(m\cdot b)$, for all $b\in B$, $m\in M$.
 By \cite[Theorem~3.15]{MenNas:ind} one thus deduces that $A\subset
X$ is Frobenius if and only if ${\rm Hom}(B,-)$ is also the left adjoint
of $\CF$.

(2) $\Leftrightarrow$ (3). By assumption,  $B$ is a finitely generated 
projective
$k$-module, $A$ is a faithfully flat $k$-module, and $X\cong B\otimes A$
as a right $A$ module. Thus $X$ is a finitely
generated projective right $A$-module (cf. \cite[I.3.6
Proposition~12]{Bou:com}). Therefore $A\subset X$ is a Frobenius
extension if and only if $X\cong {\rm Hom}_A(X,A)$ as
$(A,X)$-bimodules. One easily checks that the map $\eta_{B,A}: 
{\rm Hom}(B,A)\to {\rm Hom}_A(X, A)$ given by 
$\eta_{B,A}(f)(b\otimes a) = f(b)a$ for all $a\in A, b\in B$ is an
isomorphism of $(A,X)$-bimodules. Its
inverse is $\eta_{B,A}^{-1}(g)(b) = g(b\otimes 1_A)$, for all $b\in B$.
Combining Frobenius isomorphism with $\eta_{B,A}^{-1}$ one obtains the
required isomorphism of $(A,X)$-modules.

(3) $\Leftrightarrow$ (4). We will show that the map $\theta: {\rm
Int}(B\#_{\bar\psi}A) \to {}_A{\rm Hom}_X(X,{\rm Hom}(B,A))$, given by
$\theta(\lambda)(x) = \lambda\cdot x$, for all $x\in X$,
$\lambda\in{\rm Int}(B\#_{\bar\psi}A)$ is well-defined and bijective. For any
integral $\lambda$, the map $\theta(\lambda)$ is clearly a right
$X$-module map. Now take any $x\in X$, $a\in A$ and compute
\begin{eqnarray*}
(a\cdot\theta(\lambda))(x) & = & a\cdot\theta(\lambda)(x) =
a\cdot\lambda\cdot x \\
& = & (\lambda\cdot a)\cdot x \qquad \mbox{\rm ($\lambda$ is an
integral)}\\
& = & \theta(\lambda)(ax).
\end{eqnarray*}
Consider $\tilde{\theta} : {}_A{\rm Hom}_X(X,{\rm
Hom}(B,A))\to {\rm Int}(B\#_{\bar\psi}A)$, $\tilde{\theta}: \phi \mapsto 
\phi(1_X)$.  This map is well-defined since $\phi$ is an
$(A,A)$-bimodule map and thus for all $a\in A$,
$a\cdot \phi(1_X) = \phi(a) = \phi(1_X)\cdot a$. An easy
calculation shows that $\tilde{\theta}$ is the inverse of $\theta$.
\end{proof}

If $(A,C)_\psi$ if a factorisable entwining structure then
every $M\in
\M^C_A(\psi)$ can be viewed as an object in $\M_X$, where 
$X=C^{*op}\#_{\bar{\psi}}A$,  via the action $m\cdot
(\xi\otimes a) = m\sw 0\cdot a\<m\sw 1 ,\xi\>$, for all $m\in M$,
$a\in A$, $\xi\in C^*$. If $C$ is a finitely generated projective
$k$-module, then every
right $X$-module $M$ is an object in $\M^C_A(\psi)$ with the natural
right $A$-action and the coaction given by $m\sw 0\<m\sw 1,\xi\>
= m\cdot \xi$ for all $m\in M$, $\xi\in C^{*op}$. Therefore, if $C$ is a finitely
generated projective $k$-module then $\M_A^C(\psi)\cong \M_X$ (cf.
\cite[(1.3)]{Doi:uni}). 
Using this
category isomorphism one can reformulate Proposition~\ref{prop.frobenius1}.
Firstly, however, we need the following
\begin{definition}
Let $(A,C)_\psi$ be an entwining structure. An {\em integral} in $(A,C)_\psi$ is an
element $x=\sum_{i=1}^n a_i\otimes c_i\in A\otimes C$ such that for all
$a\in A$, $a\cdot x = x\cdot a$. The left action is the obvious one,
while the right action is as in 
Example~\ref{ex.psi.modules}(1). Explicitly, we require $\sum_{i=1}^n
aa_i\otimes c_i = \sum_{i=1}^na_i\psi(c_i\otimes a)$.
\end{definition}
If $(A,C)_\psi$ is an entwining structure of Example~\ref{Doi-Hopf}, 
then $x$ is an integral in $(A,C)_\psi$ iff it is an $H$-integral in the
sense of \cite[Definition~2.1]{CaeMil:Doi}.
\begin{example}
Let $(A,C)_\psi$ be the canonical entwining structure associated to a
coalgebra-Galois extension $A(B)^C$ as in Example~\ref{can.ex}. Then 
$x=\sum_{i=1}^n a_i\otimes c_i$ is an integral in $(A,C)_\psi$ iff for all
$a\in A$, $a\cdot
x^\tau = x^\tau\cdot a$, where $x^\tau = \sum_{i=1}^n
a_ican^{-1}(1\otimes c_i) \in A\otimes _BA$ ($A\otimes_B A$ has the
obvious $(A,A)$-bimodule structure).
\label{lemma.int}
\end{example}

The following proposition is a general $(A,C)_\psi$-module version of 
\cite[Theorem~2.4]{CaeMil:Doi}
\begin{proposition}
Let $(A,C)_\psi$ be an entwining structure and let $B=C^{*op}$. If $A$ is a faithfully
flat $k$-module and $C$ is a finitely generated projective $k$-module
then the following are equivalent:

(1) The functor $-\otimes C:\M_A\to \M_A^C(\psi)$ is the left adjoint
of the forgetful functor $\M_A^C(\psi)\to \M_A$.

(2) The extension $A\subset B\#_{\bar{\psi}} A$ is Frobenius.

(3) $C^*\otimes A\cong A\otimes C$ as $(A,B\#_{\bar{\psi}}A)$-bimodules. 

(4) $C^*\otimes A \cong A\otimes C$ as $(A,A)$-bimodules and right
$C$-comodules.

(5) There exists an integral $x=\sum_{i=1}^n a_i\otimes c_i$ in
$(A,C)_\psi$, such that the
map $\phi: C^*\otimes A\to A\otimes C$, $\xi\otimes a\mapsto
\sum_{i=1}^n a_i\psi(\xi\triangleright c_i\otimes a)$ is bijective.
\label{prop.frobenius2}
\end{proposition}

Since $\M^C_A(\psi)\cong \M_{B\#_{\bar{\psi}} A}$ 
and ${\rm Hom}(C^*, -)\cong
-\otimes C$ (for $C$ is a finitely generated projective $k$-module), the
first statement is the same as Proposition~\ref{prop.frobenius1}(1).
Since $C^*$ is a right $C$-comodule via coproduct in $C$, both 
$C^*\otimes A, A\otimes C$ are objects in $\M_A^C(\psi)$ by
Example~\ref{ex.psi.modules}. They are 
also right $B\#_{\bar{\psi}} A$ modules, via the category
isomorphism described above. $C^*\otimes A$ is
a left $A$-module via the multiplication in $B\#_{\bar{\psi}} A$, and
$A\otimes C$ is a left $A$-module via $\mu\otimes C$. Thus the statements
(3) and (4) make sense. The natural isomorphism ${\rm
Hom}(A,C^*) \cong A\otimes C$ also implies that both these
statements are equivalent to Proposition~\ref{prop.frobenius1}(3). For
the same reason the
$k$-module 
${\rm Int}(B\#_{\bar\psi}A)$ is isomorphic to the $k$-module of integrals in $(A,C)_\psi$ and
thus the statement (5) is just a reformulation of 
Proposition~\ref{prop.frobenius1}(4).

\begin{remark}
If $(A,C)_\psi$ is
an entwining structure, $A$ is a faithfully flat $k$-module and $C$ is a
projective $k$-module and if $-\otimes C$ is the
left adjoint of the forgetful functor $\M^C_A(\psi)\to \M_A$, then
$C$ is finitely generated. This follows from the fact that, in this
case,  for all $M\in \M_A$, $N\in \M_A^C(\psi)$, $m\in M$, and  
$f\in {\rm Hom}^C_A(M\otimes C,N)$,
$$
\eta_{M,N}(f)(m)=\sum_{i=1}^nf(m\cdot a_i\otimes c_i),
$$
where $\eta_{M,N}: {\rm Hom}_A^C(M\otimes C,N)\to {\rm Hom}_A(M,N)$ is a
natural isomorphism and $\sum_{i=1}^n a_i\otimes c_i = \eta_{A,A\otimes
C}(A\otimes C)(1_A)$. This can be proven following the same arguments as
in \cite[Lemma~2.3, Theorem~2.4]{CaeMil:Doi}. Therefore, 
in the case of an entwining structure associated to a
Doi-Hopf datum $(A,C,H)$ (cf. Example~\ref{Doi-Hopf}),
Proposition~\ref{prop.frobenius2} is equivalent (up to left-right
conventions, and to the assumption that $H$ has 
a bijective antipode) to
\cite[Theorem~2.4]{CaeMil:Doi}. 
\end{remark}

The following proposition generalises
\cite[Theorem~3.4]{CaeMil:Doi} to any entwining structure. 
\begin{proposition}
Let $(A,C)_\psi$ be an entwining structure 
and assume that $C$ is a projective
$k$-module. Then the following statements are equivalent:

(1) There exists $e\in C$ such that $C^*\la e = C$ and for
all $a\in A$, $\psi(a\otimes e) = e\otimes a$.

(2) $C$ is a finitely generated $k$-module and there exists a right
$C$-comodule isomorphism $\phi:C^*\to C$ such that the following
diagram
$$
\begin{CD}
A\otimes C^* @>\bar{\psi}>> C^*\otimes A \\
@VV{A\otimes\phi}V          @VV{\phi\otimes A}V \\
A\otimes C   @<{\psi}<<     C\otimes A
\end{CD}
$$
commutes ($C^*$ is a right $C$-comodule via $\xi\sw 0\<\xi\sw 1,\xi'\>
= \xi'\xi$, $\forall \xi,\xi'\in C^*$).

Furthermore, if $k$ is a field then (1) and (2) are equivalent to:

(3) $C$ is finite dimensional and there exists a non-degenerate, 
bilinear and associative form $[-,-]: C^*\otimes C^*\to k$ such that
the following diagram
$$
\begin{CD}
A\otimes C^*\otimes C^*  @>\bar{\psi}\otimes C^*>> C^*\otimes A\otimes
C^* \\
@VV{A\otimes [-,-]}V          @VV{C^*\otimes \bar{\psi}}V \\
A @<{[-,-]\otimes A}<<      C^*\otimes C^*\otimes A
\end{CD}
$$
commutes.

(4) There exists $e\in C$ such
that $e\ra C^* = C$ and for
all $a\in A$, $\psi(a\otimes e) = e\otimes a$.
\label{prop.frobenius3}
\end{proposition}
\begin{proof}
(1)$\Rightarrow$ (2).Let
$\Delta(e) = \sum_{i=1}^n c_i\otimes c'_i$. Then $C=C^*\la e\subseteq
\sum_{i=1}^n kc_i$ and thus is finitely generated. Since $C$ is also
projective, its dual $C^*$ is a finitely generated projective
$k$-module.  Consider $\phi: C^*\to C$, $\xi\mapsto \xi\la e =
e\sw 1\< e\sw 2,\xi\>$. Clearly $\phi$ is a right $C$-comodule map.  
By assumption, $\phi$ is surjective and since both
$C$ and $C^*$ are  finitely generated projective $k$-modules of the
same rank, we conclude that
$\phi$ is an isomorphism. 

Using the notation for $\bar\psi$ introduced in Section~2, we have for
all $a\in A$, $\xi\in C^*$:
\begin{eqnarray*}
\psi\circ (\phi\otimes A)\circ\bar{\psi}(a\otimes \xi) &=&
\<e\sw 2,\xi_i\>\psi(e\sw 1\otimes a^i) \\
& = & \<e\sw 2^\alpha,\xi\>\psi(e\sw 1\otimes a_\alpha)
\qquad \mbox{\rm (by the definition of $\bar{\psi}$)}\\
& = & \<e^\alpha\sw 2,\xi\>a_\alpha\otimes e^\alpha\sw 1
\qquad\quad \mbox{\rm (by Definition~\ref{ent})}\\
& = & a\otimes e\sw 1\<e\sw 2,\xi\> \qquad \qquad \quad\mbox{\rm (for
$\psi(e\otimes a) = a\otimes e$)}\\
& = & a\otimes \phi(\xi).
\end{eqnarray*}

(2) $\Rightarrow$ (1). Given $\phi:C^*\to C$, define $e = \phi(\eps)$.
Since $\phi$ is a right $C$-comodule map 
and $\eps$ is the unit in $C^*$ one easily finds 
that for all $\xi\in C$, $\phi(\xi) = \xi\la
e$. Thus $C^*\la e = C$. The
commutativity of the diagram in (2) together with the fact that for all
$a\in A$, $\bar{\psi}(a\otimes \eps) = \eps\otimes a$ now imply that
$\psi(e\otimes a) = a\otimes e$, as required.

(2) $\Rightarrow$ (3). Given $\phi: C^*\to C$, define $[-,-]: C^*\otimes
C^*\to k$, by $[-,-]: \xi\otimes \xi'\mapsto \<\phi(\xi'),\xi\>$. Since
$\phi$ is right $C$-colinear, it is left $C^*$-linear and thus for all
$\xi, \xi',\xi''\in C^*$
$$
[\xi\xi',\xi''] = \<\phi(\xi''),\xi\xi'\> = \<\xi'\la\phi(\xi''),\xi\>
= \<\phi(\xi'\xi''),\xi\> = [\xi,\xi'\xi''],
$$
so that $[-,-]$ is associative. Since $\phi$ is onto, $C^*\ni\xi =0$ if
and only if for all $\xi'\in C^*$, $[\xi,\xi'] = \<\phi(\xi'),\xi\> =
0$. Assume now that there exists $\xi\in C^*$ such that for all $\xi'\in
C^*$ we have $[\xi',\xi] =0$. Let $\{\xi_n,c_n\}_{n=1}^N$ be a  
dual basis in
$C$. Then $0=\sum_{n=1}^N[\xi_n,\xi]c_n = \sum_{n=1}^N\<\phi(\xi),
\xi_n\>c_n =
\phi(\xi)$, and thus $\xi =0$ since $\phi$ is bijective. Therefore the
form $[-,-]$ is non-degenerate. Finally, take any $\xi,\xi'\in C^*$,
$a\in A$ and, using the same notation for $\bar{\psi}$ as
before and  the commutativity of the diagram in statement (2), compute
$
\<\phi(\xi'_i),\xi_j\> a^{ji}=
\<\phi(\xi'_i)^\alpha,\xi\>a_\alpha^i = \<\phi(\xi'),\xi\>
$. This
 implies the commutativity of the diagram in (3).

Now assume that $k$ is a field.

(3) $\Rightarrow$ (1). Given a non-degenerate, associative form $[-,-]$
on $C^*$, and a dual basis $\{\xi_n,c_n\}_{n=1}^N$ in $C$ define
$e=\sum_{n=1}^N [\eps, \xi_n]c_n$. Then for all
$\xi,\xi'\in C^*$ one has $[\xi,\xi'] = \<e,\xi\xi'\> = \<\xi'\la
e,\xi\>$. The non-degeneracy of $[-,-]$ implies that the map $C^*\to C$,
$\xi\mapsto \xi\la e$ is injective. Since both $C$ and $C^*$ have the same
dimension this map is also surjective. Now notice that for all $a\in A$,
$\sum_{n=1}^N \psi(c_n\otimes a)\otimes \xi_n =
\sum_{n=1}^N a^i\otimes c_n\otimes \xi_{ni}$, where we use the same
notation for $\bar{\psi}$ as before. Using this fact, as well
as $\bar{\psi}(\eps\otimes a) = a\otimes \eps$ and the commutativity of
the diagram in (3) one computes
$$
\psi(e\otimes a) = \sum_{n=1}^N [\eps,\xi_n]\psi(c_n\otimes a) =
\sum_{n=1}^N[\eps_i,\xi_{nj}]a^{ij}\otimes c_n = \sum_{n=1}^N a\otimes
[\eps,\xi_n]c_n = a\otimes e.
$$

(4) $\Rightarrow$ (1). The similar argument as in the proof of the first
implication shows that $C$ is finite dimensional. Since $e\ra C^* =C$,
we have for all $\xi\in C$,
$
\<C,\xi\> = \<e\ra C^*,\xi\> = \<\xi\la e,C^*\>.
$
Therefore, $\xi\la e =0$ implies $\xi =0$ so that the map $\xi\mapsto
\xi\la e$ is injective. Since both $C$ and $C^*$ are of the same finite
dimension, we conclude that $C^*\la e = C$. Using similar argument one
shows the implication (1) $\Rightarrow$ (4).
\end{proof}

Statement (1) of Proposition~\ref{prop.frobenius3} takes a very
simple form in the case of the canonical  entwining structure
$(A,C)_\psi$ 
associated to a $C$-Galois extension $A(B)^C$. In this case 
$e\in C$ has the property that
for all $a\in A$, $\psi(e\otimes a) = a\otimes e$, if and only if for
all $a\in A$, $a\tau(e) = \tau(e) a$, where $\tau(c) =
\can^{-1}(1\otimes c)$.
This is a 
simple consequence of the definition of
the canonical entwining structure in Example~\ref{can.ex}.

\section{Maschke-type theorems}
The classical Maschke's Theorem states that a group ring of a finite
group is semisimple if and only if the characteristic of the field does
not divide the order of the group. Several generalisations of Maschke's
theorem to Hopf algebras and comodule algebras are known
\cite{LarSwe:ass} \cite{Doi:Mas} \cite{BlaMon:cro}.  In 
\cite{CaeMil:Mas} a Maschke-type
theorem was formulated for Doi-Hopf modules.
 Following \cite{CaeMil:Mas} we define
\begin{definition}
Let $(A,C)_\psi$ be an entwining structure. A $k$-module map 
$\phi: C\to C^*\otimes A$ is called a {\em normalised integral map} in
$(A,C)_\psi$ if  

(1) $\mu\circ(A\otimes \ev_C\otimes A)\circ (\psi\otimes
\phi)\circ(C\otimes \psi) = (\ev_C\otimes
\mu)\circ(C\otimes\phi\otimes A)$

(2) $(\ev_C\otimes A\otimes C)\circ(C\otimes\phi\otimes C)\circ
(C\otimes\Delta) = \psi\circ(C\otimes \ev_C\otimes A)\circ
(\Delta\otimes \phi)$,

(3) 
$(\ev_C\otimes A)\circ(C\otimes\phi)\circ\Delta = 1_A\circ\eps$.
\label{int.map}
\end{definition}
If $(A,C)_\psi$ is a factorisable entwining structure then conditions
(1) and (2) in Definition~\ref{int.map} state that $\phi$ is an
$(A,C^{*op}\#_{\bar\psi}A)$-bimodule map. Condition (3) is a
normalisation condition. In particular, for an entwining structure of
Example~\ref{Doi-Hopf}, Definition~\ref{int.map} is equivalent to
\cite[Definition~2.1]{CaeMil:Mas} (up to the left-right conventions and
the bijectivity of an antipode). The following generalises
\cite[Theorem~2.5]{CaeMil:Mas}
\begin{theorem}
If there exists a normalised integral map in $(A,C)_\psi$ then a 
morphism in
$\M_A^C(\psi)$ which has a section (resp.\ retraction) in $\M_A$, has a
section (resp.\ retraction) in $\M_A^C(\psi)$.
\label{thm.maschke}
\end{theorem}
\begin{lemma}
Let $M,N\in \M_A^C(\psi)$, $g\in {\rm Hom}_A(M,N)$ and let $\phi$ be a
normalised integral map in $(A,C)_\psi$. Then 
$\tilde{g} :M\to N$, $\tilde{g} = \rho_N\circ(N\otimes \ev_C\otimes
A)\circ(\rho^N\circ g\otimes \phi)\circ \rho^M$
 is a morphism in $\M_A^C(\psi)$.
\label{lemma2}
\end{lemma}
\begin{proof}
We use the Sweedler-like notation $\phi(c) = c\su 1\otimes c \su 2\in
C^*\otimes A$ (summation understood), for any
$c\in C$. Then $\tilde{g}(m) =  
g(m\sw 0)\sw 0\cdot m\sw 1\su 2\<g(m\sw 0)\sw 1, m\sw
1\su 1\>$, for all $m\in M$. 
To show that $\tilde{g}$ is right $A$-linear, take any $m\in M$,
$a\in A$ and compute
\begin{eqnarray*}
\tilde{g}(m\cdot a) & = & g((m\cdot a)\sw 0)\sw 0\cdot 
(m\cdot a)\sw 1\su 2 \<g((m\cdot a)\sw 0)\sw 1, (m\cdot a)\sw 1\su 1\> \\
& = & g(m\sw 0\cdot a_\alpha)\sw 0\cdot m\sw
1^\alpha\su 2\< g(m\sw 0\cdot a_\alpha)\sw 1, m\sw 1^\alpha\su 1\>
\;\qquad \mbox{\rm ($M\in \M_A^C(\psi)$)}\\
& = & (g(m\sw 0)\cdot a_\alpha)\sw 0\cdot m\sw
1^\alpha\su 2\< (g(m\sw 0)\cdot a_\alpha)\sw 1, m\sw 1^\alpha\su 1\>
\quad \mbox{\rm ($g$ is $A$-linear)}\\
& = & g(m\sw 0)\sw 0\cdot a_{\alpha\beta}m\sw
1^\alpha\su 2\< g(m\sw 0)\sw 1^\beta, m\sw 1^\alpha\su 1\> 
\qquad\qquad\mbox{\rm ($N\in \M_A^C(\psi)$)}\\
& = & g(m\sw 0)\sw 0\cdot m\sw 1\su 2a \<g(m\sw 0)\sw 1, m\sw 1\su 1\>
\quad\qquad\qquad
\qquad\mbox{\rm (Def.~\ref{int.map}(1))}\\
& = & \tilde{g}(m)\cdot a.
\end{eqnarray*}
Next, for all $m\in M$ we have,
\begin{eqnarray*}
\rho^N(\tilde{g}(m)) & = & (g(m\sw 0)\cdot m\sw 1\su 2)\sw 0 \<g(m\sw
0)\sw 1, m\sw 1\su 1\> \otimes (g(m\sw 0)\cdot m\sw 1\su 2)\sw 1 \\
& = & g(m\sw 0)\sw 0\cdot m\sw 1\su 2_\alpha \<g(m\sw
0)\sw 2, m\sw 1\su 1\> \otimes g(m\sw 0)\sw 1^\alpha\\
& = & g(m\sw 0)\sw 0\cdot m\sw 1\su 2 \<g(m\sw
0)\sw 1, m\sw 1\su 1\> \otimes m\sw 2 \qquad \qquad
\mbox{\rm (Def.~\ref{int.map}(2))}\\
& = & \tilde{g}(m\sw 0)\otimes m\sw 1,
\end{eqnarray*}
where we used that $N\in \M_A^C(\psi)$ to derive the second equality.
\end{proof}

{\sl Proof of Theorem~\ref{thm.maschke}.}
Let $M,N\in \M_A^C(\psi)$, and assume that $f\in {\rm Hom}_A^C(N, M)$
has a section $g\in {\rm Hom}_A(M, N)$. Let
$\tilde{g}\in{\rm Hom}^C_A(M,N)$ be as in Lemma~\ref{lemma2}. Then for
all $m\in M$, 
\begin{eqnarray*}
f\circ \tilde{g}(m) & = & f(g(m\sw 0)\sw 0\cdot m\sw 1\su 2) \<g(m\sw
0)\sw 1, m\sw 1\su 1\>\\
& = & f(g(m\sw 0))\sw 0 m\sw 1\su 2 \<f(g(m\sw
0))\sw 1, m\sw 1\su 1\>\quad
\mbox{\rm ($f\in {\rm Hom}_A^C(N, M)$)}\\
& = & m\sw 0\cdot m\sw 2\su 2\<m\sw 1,m\sw 2\su 2\> = m \qquad
\mbox{\rm ($g$ is a section of $f$, Def.~\ref{int.map}(3))}
\end{eqnarray*}
Similar computation shows that if $g$ is a retraction of $f$ then
so is $\tilde{g}$. 
\endproof

\begin{corollary}
If there is a normalised integral map in $(A,C)_\psi$, then 

(1) Every 
object in $\M_A^C(\psi)$ which is semisimple as an object in $\M_A$ is
semisimple as an object in $\M_A^C(\psi)$. 

(2) Every object in
$\M_A^C(\psi)$ which is projective (resp.\ injective) as a 
right $A$-module is
a projective (resp.\ injective) object in $\M_A^C(\psi)$. 

(3) If $C$ is a
flat $k$-module then $M\in \M_A^C(\psi)$ is
projective as a right $A$-module if and only if there exists $V\in \M^C$
such that $M$ is a direct summand of $V\otimes A$ in $\M_A^C(\psi)$
($V\otimes A$ is an entwined module by Example~\ref{ex.psi.modules}(2)).
\end{corollary}
\begin{proof}
The first statement follows immediately from Theorem~\ref{thm.maschke}.
Assertions (2) and (3) can be proven by the same
method as \cite[Corollary~2.9]{CaeMil:Mas}.
In particular to prove (2) it is useful to observe that if $f$ is a
morphism in $\M^C_A(\psi)$, then $\tilde{f} = f$, where $\tilde{f}$ is
constructed in Lemma~\ref{lemma2}. To prove (3) one uses 
that the forgetful functor 
$\M_A^C(\psi)\to \M_A$ is the left adjoint of the exact 
functor $-\otimes C:
\M_A \to \M_A^C(\psi)$ (cf. \cite[Section~3]{Brz:mod}).
\end{proof}

We can dualise the above construction to derive the dual version of
Theorem~\ref{thm.maschke}. 
\begin{definition}
Let $(A,C)_\psi$ be an entwining structure with $A$ finitely generated
projective $k$-module. Let $\coev_A :k\to A\otimes A^*$ be a coevaluation
map, i.e., $\coev_A :\kappa \mapsto\kappa\sum_{i\in I}a_i\otimes a_i^*$,
where $\{a_i,a^*_i\}_{i\in I}$ is a dual basis in $A$.
A $k$-module map 
$\phi: A^*\otimes C\to A$ is called a {\em normalised cointegral map} in
$(A,C)_\psi$ if  

(1) $(A\otimes\psi)\circ(\psi\otimes\phi)\circ(C\otimes\coev_A\otimes
C)\circ\Delta = (A\otimes\phi\otimes C)\circ(\coev_A\otimes
\Delta)$,
 
(2) $(A\otimes\mu)\circ(A\otimes\phi\otimes A)\circ(\coev_A\otimes
C\otimes A) = (\mu\otimes \phi)\circ(A\otimes\coev_A\otimes
C)\circ\psi$, 

(3) 
$\mu\circ(A\otimes \phi)\circ(\coev_A\otimes C)= 1_A\circ\eps$.
\label{coint.map}
\end{definition}

Conditions (1) and (2) in Definition~\ref{coint.map} can be understood
as follows. $A^*$ is a right $A$-module with 
the action $\<a',a^*\cdot a \>= \<aa',a^*\>$, for
any $a,a'\in A$, $a^*\in A^*$. Thus $A^*\otimes C$ is an
entwined module as in Example~\ref{ex.psi.modules}(1). 
Furthermore, $\psi$ induces the map 
$\hat{\psi} : A^*\otimes
C\to C\otimes A^*$,  $a^*\otimes c\mapsto \sum_{i\in I} 
c^\alpha 
\otimes a^*_i \<a_{i\alpha}, a^*\>$.
  Using this map one defines a
left coaction of $C$ on $A^*\otimes C$ as $(\hat{\psi}\otimes
C)\circ (A^*\otimes\Delta)$. The $k$-module $C\otimes A$ has a
right $A$ module and right $C$-comodule structure as in
Example~\ref{ex.psi.modules}(2), and a left $C$-comodule structure,
$\Delta\otimes C$. Conditions (1) and (2) of Definition~\ref{coint.map} 
are equivalent to the existence of
$\tilde{\phi}\in {}^C{\rm Hom}^C_A( A^*\otimes C , C\otimes A)$ (then 
$\phi = (\eps\otimes A)\circ\tilde{\phi}$).

\begin{theorem}
If there exists a normalised cointegral map in $(A,C)_\psi$,
then any morphism in
$\M_A^C(\psi)$ which  has a section (resp.\ retraction) in $\M^C$, has a
section (resp.\ retraction) in $\M_A^C(\psi)$.
\label{thm.maschke*}
\end{theorem}
\begin{lemma}
Let $M,N\in \M_A^C(\psi)$, $g\in {\rm Hom}^C(M,N)$ and let $\phi$ be a
normalised cointegral map in $(A,C)_\psi$. Then 
$\tilde{g} :M\to N$, $\tilde{g} =
\rho_N\circ(g\circ\rho_M\otimes\phi)\circ(M\otimes\coev_A\otimes C)\circ
\rho^M$ 
is a morphism in $\M_A^C(\psi)$.
\label{lemma2*}
\end{lemma}
\begin{proof}
Dual to the proof of Lemma~\ref{lemma2}.
\end{proof}

{\sl Proof of Theorem~\ref{thm.maschke*}.}
Let $M,N\in \M_A^C(\psi)$, and assume that $f\in {\rm Hom}_A^C(N, M)$
has a section $g\in {\rm Hom}^C(M, N)$. Let $\tilde{g}\in {\rm Hom}_A^C(N,
M)$ be as in Lemma~\ref{lemma2*}. Explicitly, for all $m\in M$, $\tilde{g}(m)= \sum_{i\in I} g(m\sw
0\cdot a_i)\cdot 
\phi(a^*_i\otimes m\sw 1)$, where $\{a_i,a^*_i\}_{i\in I}$ is a dual
basis in $A$. Thus
\begin{eqnarray*}
f\circ \tilde{g}(m) & = & \sum_{i\in I} f(g(m\sw 0\cdot a_i)\cdot
\phi(a^*_i \otimes m\sw 1)) \\
& = & \sum_{i\in I}f(g(m\sw 0\cdot a_i))\cdot
\phi(a^*_i \otimes m\sw 1) \qquad
\mbox{\rm ($f$ is right $A$-linear)}\\
& = & \sum_{i\in I} m\sw 0\cdot a_i\phi(a^*_i\otimes m\sw 1) = m \qquad
\mbox{\rm ($g$ is a section of $f$, (Def.~\ref{coint.map}))}\\
\end{eqnarray*}
Similar computation shows that if $g$ is a retraction of $f$ then
so is $\tilde{g}$. 
\endproof

Theorem~\ref{thm.maschke*} implies that if there is a 
normalised cointegral map in $(A,C)_\psi$, then every 
object in $\M_A^C(\psi)$ which is semisimple as an object in $\M^C$ is
semisimple as an object in $\M_A^C(\psi)$. Furthermore every object in
$\M_A^C(\psi)$ which is a projective (resp.\ injective) object in 
$\M^C$ is
projective (resp.\ injective) in $\M_A^C(\psi)$. 
\note{Also, 
$M\in \M_A^C(\psi)$ is
an injective object in $\M^C$ if and only if there exists $V\in \M_A$
such that $M$ is a direct summand of $V\otimes C$ in $\M_A^C(\psi)$
($V\otimes C$ is an entwined module by Example~\ref{ex.psi.modules}(1)).
The latter follows from the fact that the forgetful functor
$\M_A^C(\psi)\to \M^C$ is the right adjoint of the exact 
functor $-\otimes A:
\M^C \to \M_A^C(\psi)$ (cf. \cite[Section~3]{Brz:mod}).}

\end{document}